\input amstex

\documentstyle{amsppt}
\nologo
\NoBlackBoxes

\magnification\magstep1
\pagewidth{32.622pc}
\pageheight{44.8pc}

\topmatter

\title
A simple proof of Bernstein-Lunts equivalence
\endtitle
\author 
Pavle Pand\v zi\'c
\endauthor 
\thanks Pavle Pand\v zi\'c, 
Department of Mathematics,
University of Zagreb,
PP 335, 10002 Zagreb, Croatia; email adress: pandzic\@math.hr
\endthanks
\subjclass \nofrills 2000 {\it Mathematics Subject Classification}. 
Primary 22E46
\endsubjclass
\abstract We give an easy proof of the Bernstein-Lunts equivalence of ordinary
and equivariant derived categories of Harish-Chandra modules. This proof
requires no boundedness assumptions. In the appendix
we collect some needed, but not completely standard facts from homological
 algebra.
\endabstract
\endtopmatter

\document

\def\Hom{\operatorname{Hom}}

\def\Int{\operatorname{Int}}
\def\Ind{\operatorname{Ind}}

\def\pro{\operatorname{pro}}

\def\ad{\operatorname{ad}}

\def\Ob{\operatorname{Ob}}
\def\For{\operatorname{For}}

\def\deg{\operatorname{deg}}

\def\A{\Cal A}
\def\B{\Cal B}
\def\C{\Cal C}
\def\D{\Cal D}
\def\E{\Cal E}
\def\U{\Cal U}
\def\M{\Cal M}
\def\Sl{\Cal S}
\def\T{\Cal T}
\def\N{\Cal N}

\def\L{\Cal L}
\def\R{\Cal R}

\def\g{\frak g}
\def\k{\frak k}

\def\Xd{X^\cdot}

\def\Vd{V^\cdot}

\def\Cd{C^\cdot}

\def\dotimes{{\overset L \to \otimes}}

\def\eps{\varepsilon}
\def\no{\noindent}

\def\pf{\demo{Proof}}
\head Introduction
\endhead

The aim of this paper is to give an easy proof of Bernstein-Lunts
equivalence, the main result of \cite{BL2}. This easy proof has an 
additional advantage: it goes through without any boundedness assumptions.

To explain the setting, let $\g$ be a complex Lie algebra, and let $K$
be a complex algebraic group acting on $\g$ via a morphism 
$\phi:K@>>>\Int(\g)$, so that the differential of $\phi$ defines an 
embedding of $\k$ into $\g$. Then $(\g,K)$ is called a Harish-Chandra
pair. Let $\M(\g,K)$ be the category of Harish-Chandra modules for the
pair $(\g,K)$; these are modules simultaneously for $\g$ and $K$, with
the usual compatibility conditions. Let $D(\M(\g,K))$ be the 
derived category of the abelian category $\M(\g,K)$. There is another 
related notion, the equivariant derived category 
$D(\g,K)$ for the pair $(\g,K)$. Objects of this category are equivariant
$(\g,K)$-complexes; these are complexes of ``weak'' Harish-Chandra modules
(with weakened compatibility conditions), 
endowed with a family of homotopies $i_\xi$, $\xi\in\k$ satisfying
certain properties; for the precise definition, see \cite{BB}, \cite{G}, 
\cite{BL2},
\cite{MP} or \cite{P2}. Equivalently, equivariant $(\g,K)$-complexes are
$(\g,K,N(\k))$-modules where the differential graded (DG) algebra $N(\k)$ 
is the standard complex of $\k$. This setting is explained in \S 1 below.

There is an obvious functor $D(\M(\g,K))@>>>D(\g,K)$: a complex of
Harish-Chandra modules can be viewed as an equivariant $(\g,K)$-complex
with all $i_\xi=0$. One of the main results of \cite{BL2}, Theorem 1.10,
asserts that $Q$ is an equivalence of bounded derived categories.

The proof in \cite{BL2} is rather complicated; it uses K-injective 
resolutions and some dualizing arguments. We first give another proof
for reductive $K$, using K-projectives instead of K-injectives. This
proof is similar to the proof of an analogous result for DG modules
over DG algebras in \cite{BL1}, 10.12.5.1. No boundedness needs to be
assumed; see the end of \S 1 for an explanation why.
For non-reductive $K$, the result now follows by using the arguments
of \cite{MP, \S 2}.

In the following let us briefly describe the contents of the paper.
In \S 1 we describe a construction of enough K-projective equivariant
$(\g,K)$-complexes. In \S 2 we give the proof of Bernstein-Lunts
equivalence (Theorem 2.4 and Corollaries 2.5 and 2.6).

In the appendix we collect some facts about homological algebra needed
in the paper. These are known, but they are not readily available in
the literature, at least not in the form needed.

In \S A1 we explain the definition and construction of derived functors
in the setting of a triangulated category and its localization, which
are not necessarily the homotopic and derived categories of an abelian
category. This is necessary to study equivariant derived categories.
In \S A2 we collect some facts about adjoint functors, in particular the
ones related to homological algebra. These are used over and over, both
in this paper and in \cite{P2}.

\head 1. K-projectives in equivariant derived categories
\endhead

Let us start be recalling the definition of $(\A,K,\D)$-modules from \cite{P2}.
Let $K$ be a complex algebraic group with Lie algebra $\k$. Let $\A$ be an 
associative algebra over $\Bbb C$, with an algebraic action $\phi$ of $K$, and 
a $K$-equivariant Lie algebra morphism $\psi:\k@>>>\A$, such that the 
differential of $\phi$ satisfies
$$
d\phi(\xi)(a)=[\psi(\xi),a],\quad \xi\in\k, a\in\A.
$$
Let $\D$ be a DG algebra over $\Bbb C$ with an algebraic action $\chi$ of $K$
and a morphism $\rho:\k@>>>\D$ of DG Lie algebras, satisfying analogous
conditions. 

An $(\A,K,\D)$-module is a complex $V$ of vector spaces, with an action $\pi$
of $\A$ by chain maps, an algebraic action $\nu$ of $K$ by chain maps, and 
a DG action $\omega$ of $\D$, such that $\pi$ and $\omega$ commute and are both 
$K$-equivariant, and such that $\pi+\omega=\nu$ on $\k$.

We denote the abelian category of $(\A,K,\D)$-modules by $\M(\A,K,\D)$; the 
morphisms are chain maps which preserve all the actions.
One defines the homotopic category $K(\A,K,\D)$ as in \cite{P2}, Section 2.3;
it is a triangulated category. Localizing with respect to quasiisomorphisms,
we get to the equivariant derived category $D(\A,K,\D)$. This definition is
due to Beilinson and Ginzburg; an analogous definition in geometric setting
is due to Bernstein and Lunts. See \cite{BB}, \cite{G}, \cite{BL1} and \cite{BL2}. 

The main example is when $\D=N(\k)$, the standard complex of $\k$.

Now we want to show that in case $K$ is reductive, the category
$K(\A,K,\D)$ has enough K-projectives (see Section A1.4). The proof is 
analogous to the proof of the same fact for the category of DG modules over a 
DG algebra $\D$ from \cite{BL1}, 10.12.2.
The main idea is familiar: one uses the fact that there are enough projectives
in the category of weak $(\A,K)$-modules, hence there are enough K-projectives
in the homotopic category of complexes $K(\M(\A,K)_w)=K(\A,K,\U(\k))$ 
(see A1.4.7). 
Now one constructs K-projectives in $K(\A,K,\D)$ applying the change
of DG algebras from \cite{P2}, Section 2.5, that is, the functor
$V\longmapsto V\otimes_{\U(\k)}\D$. 

The construction goes as follows. 
Let $V$ be an $(\A,K,\D)$-module. Forgetting the $\D$-action, we get a
complex of weak $(\A,K)$-modules. Let $Q@>s>>V$ be a K-projective resolution
of $V$ in the category $C(\M(\A,K)_w)$. We can assume that $s$ is surjective.
Let $P_0=Q\otimes_{\U(k)}\D$. Then $s\otimes 1:P_0@>>>V\otimes_{\U(k)}\D$
is still surjective since the functor $-\otimes_{\U(k)}\D$ is right exact being
a left adjoint. Furthermore, the adjunction morphism 
$\Psi_V:V\otimes_{\U(k)}\D@>>>V$, which is given by 
$$
\Psi_V(v\otimes x)=(-1)^{\deg v\deg x}\omega_V(^\iota x)v,
$$
is clearly also surjective. So the composition 
$\eps_0=\Psi_V\circ (s\otimes 1):P_0@>>>V$
is also surjective (in the category $\M(\A,K,\D)$). We claim that $\eps_0$ is 
also surjective on the level of cohomology.
Since $s:Q@>>>V$ is a quasiisomorphism, each element of $H(V)$ has a 
representative of the form $s(q)$ for some cycle $q\in Q$. However, $q\otimes 1$ 
is a cycle in $Q\otimes_{\U(k)}\D$ since $d_\D(1)=0$, and 
$$
\eps_0(q\otimes 1)=\psi_V(s(q)\otimes 1)=s(q),
$$
so the cohomology class of $s(q)$ is in the image of $H(\eps_0)$.

Let $K_0$ be the kernel of $\eps_0$. Then from the long exact sequence of 
cohomology corresponding to the short exact sequence $0@>>>K_0@>>>P_0@>>>V@>>>0$
we see that $H(K)$ is the kernel of $H(\eps_0)$.

We now repeat the above discussion for $K_0$ instead of $V$ and proceed 
inductively. In this way we get a resolution
$$
\dots@>{\eps_{-2}}>>P_{-1}@>{\eps_{-1}}>>P_0@>{\eps_0}>>V
$$
of $V$ by equivariant complexes, which induces a resolution of cohomology of $V$. 

Let us now consider the complex
$$
\dots@>{\eps_{-2}}>>P_{-1}@>{\eps_{-1}}>>P_0@>>>0\dots \tag{\dag}
$$
of equivariant complexes. We want to consider the total complex of $P^\cdot$.
The construction is as follows.

Let 
$$
\Vd=\dots@>>>V^{-1}@>{\delta_{-1}}>>V^0@>{\delta_0}>>V^1@>>>\dots
$$
be a complex of $(\A,K,\D)$-modules; in particular, we can view it as a double 
complex of weak $(\A,K)$-modules. We define the total complex $s(\Vd)$ in the 
following way. As a graded $(\A,K,\D)$-module, it is the direct sum
$$
s(\Vd)=\oplus_{i=-\infty}^\infty V^i[-i];
$$
so in particular, $s(\Vd)^k=\oplus_{i=-\infty}^\infty (V^i)^{k-i}$. The 
differential $d$ is given by
$$
d^k|_{(V^i)^{k-i}}=\delta_i\oplus d^{k-i}_{V^i[-i]}=
(-1)^i\delta_i\oplus d^{k-i}_{V^i}.
$$
This is one of the standard ways to define a differential on the total complex 
of a double complex. 
Therefore, $s(\Vd)$ as defined above is a complex and a graded $(\A,K,\D)$-module.
 It is now easy to see that the DG property is satisfied, so that $s(\Vd)$
is actualy an $(\A,K,\D)$-module.

In particular, returning to our complex $P^\cdot$, we see that its total
complex $P=s(P^\cdot)$ is an $(\A,K,\D)$-module. Clearly,
$\eps_0$ induces a morphism $\eps:P@>>>V$. We want to see that
$\eps$ is a quasiisomorphism. It is enough to check this on the level of 
complexes
of vector spaces. Then $P$ is the total complex of the double complex $(\dag)$.
Since $(\dag)$ is a left half-plane double complex, the cohomology of $P$ can be 
computed from the second spectral sequence of $(\dag)$, that is the one starting
with columns. The $E_1$ term of this spectral sequence is
$$
\dots@>>>H(P_{-2})@>>>H(P_{-1})@>>>H(P_0)@>>>0\dots
$$
This is exact except at degree 0, and there cohomology is isomorphic to $H(V)$.
So the spectral sequence degenerates at $E_2$ and the cohomology of $P$ is
isomorphic to $H(V)$. It is clear from this discussion that the isomorphism
is induced by $\eps$. So indeed $\eps$ is a quasiisomorphism. 

Finally, we want to show that $P$ is a K-projective $(\A,K,\D)$-module.
In other words, we need to prove

\proclaim{1.1. Lemma} Let $\Vd$ be a complex of $(\A,K,\D)$-modules
bounded from above, such that each $V^i$ is K-projective.
Then $s(\Vd)$ is K-projective.
\endproclaim

\pf For any complex $\Vd$ of $(\A,K,\D)$ modules, it follows from the 
definition of the differential of $s(\Vd)$ that for any $p\in\Bbb Z$,
$$
F_p s(\Vd)=\oplus_{i=p}^\infty V^i[-i]
$$
is an $(\A,K,\D)$-submodule of $s(\Vd)$. It is now clear that in this way we
get a decreasing exhaustive Hausdorff filtration of $s(\Vd)$.

Suppose now $\Vd$ is bounded above by degree $m$, and each $V^i$ is K-projective. Then $F_p s(\Vd)=0$ for $p>m$, and
by defining $\tilde F_p s(\Vd)=F_{m+1-p}s(\Vd)$ we get into the situation of
1.3 below, so $s(\Vd)$ is K-projective. Namely, the graded pieces corresponding
to the above filtration are translates of $V^i$'s.
\qed\enddemo

To finish the proof of existence of enough K-projectives, it remains to
prove that (under certain conditions) a filtered $(\A,K,\D)$-module such
that the corresponding graded pieces are K-projective, is itself K-projective.
For finite filtrations this is proved in \cite{P2}, 2.6.4; the above
filtration is however infinite.

To get the desired result, we need a technical result about cones. It is
implicit in \cite{BL1}, 10.12.2.6.

Let $f:V@>>>W$ be a morphism in $\M(\A,K,\D)$. Let $\phi:C_f@>>>Z$ be another
such morphism. Using the decomposition $C_f=T(V)\oplus W$, we can write $\phi$ as $\pmatrix \phi_1 & \phi_2\endpmatrix$, where $\phi_1:T(V)@>>>Z$ is a graded $(\A,K,\D)$-morphism of degree 0, while $\phi_2:W@>>>Z$ is a morphism in $\M(\A,K,\D)$, i.e., also a chain map.

Writing out $d_Z\phi=\phi d_f$ as matrices, we get
$$
\quad  d_Z\phi_1=\phi_1 d_{T(V)}+\phi_2T(f);   \tag{*}
$$
the other matrix entry just shows that $\phi_2$ is a chain map.

\proclaim{1.2. Lemma} Assume that $\phi_2$ is homotopic to 0 via a homotopy $h_2$, i.e., $h_2:W@>>>Z$ is a graded morphism of degree $-1$ such that
$$
\quad \phi_2=h_2 d_W+d_Z h_2.  \tag{**}
$$
Let us denote the graded morphism from $T(W)$ to $Z$ of degree 0 defined by 
$h_2$ again by $h_2$. Then:
\roster
\item"(i)" $\phi_1-h_2T(f):T(V)@>>>Z$ is a morphism in $\M(\A,K,\D)$;
\item"(ii)" If $\phi_1-h_2T(f)$ is homotopic to 0, then $h_2$ extends to a 
homotopy $h$ from $\phi$ to 0.
\endroster
\endproclaim
\pf (i) Using (*) and (**), we see
$$
\multline
d_Z(\phi_1-h_2T(f))=d_Z\phi_1-d_Z h_2 T(f)=\\
\phi_1d_{T(V)}+\phi_2 T(f)-(\phi_2 T(f)-h_2 d_W T(f))=
\phi_1 d_{T(V)}-h_2 T(f)d_{T(V)}=\\
(\phi_1-h_2T(f))d_{T(V)}.
\endmultline
$$
(ii) Let $h_1:T(V)@>>>Z$ be a graded morphism of degree $-1$ such that
$$
\quad  \phi_1-h_2T(f)=h_1d_{T(V)}+d_Zh_1,  \tag{***}
$$
i.e., $h_1$ is a homotopy from $\phi_1-h_2T(f)$ to $0$. Let $h=\pmatrix h_1 & h_2  \endpmatrix :C_f@>>>Z$. It is clearly a graded morphism of degree -1. Using
(**) and (***), we see that
$$
\phi=\pmatrix \phi_1  &  \phi_2 \endpmatrix =\pmatrix h_1 & h_2  \endpmatrix
\pmatrix d_{T(V)}  &   0  \\
          T(f)     &  d_W
\endpmatrix
+d_Z \pmatrix h_1 & h_2 \endpmatrix,
$$
which shows that $h=\pmatrix h_1 & h_2 \endpmatrix$ is a homotopy from $\phi$
to 0.
\qed\enddemo

This has the following consequence, which is exactly the property
of K-projectives we need to finish the proof of 1.1.

\proclaim{1.3. Theorem} Let $V$ be an $(\A,K,\D)$-module. Let
$$
0=F_0V\subset F_1V\subset F_2V\subset\dots
$$
be an increasing exhaustive filtration of $V$ by $(\A,K,\D)$-submodules (exhaustive means that $V=\cup_iF_iV$). Assume
that the $(\A,K,\D)$-modules
$$
Gr_iV=F_iV/F_{i-1}V,\quad i=1,2,\dots
$$
are K-projective, and that as graded modules, $F_iV\cong F_{i-1}V\oplus Gr_iV$
for any $i$. Then $V$ is K-projective.
\endproclaim
\pf Let $Z$ be an acyclic $(\A,K,\D)$-module, and $f:V@>>>Z$ a morphism in the
category $\M(\A,K,\D)$. We have to prove that $f$ is homotopic to 0.

Clearly, $V$ is the direct limit of the direct system $(F_iV)$ in $\M(\A,K,\D)$, and $f$ is the direct limit of the morphisms $f_i=f|_{F_iV}$. The same is
true in the category $\M^{GR}(\A,K,\D)$; there actually $V=\oplus_{k=1}^{\infty} Gr_kV$ and $F_iV=\oplus_{k=1}^i Gr_kV$ for $i>0$. Therefore, it is enough to construct
homotopies $h_i:F_iV@>>>Z$ from $f_i$ to 0 for every $i$, which are compatible,
i.e., $h_i|_{F_{i-1}V}=h_{i-1}$. Then they define a graded morphism $h:V@>>>Z$
of degree -1, and $h$ is a homotopy from $f$ to 0 because 
$$
f=d_Z h+hd_V
$$
follows from the fact that for every $i$
$$
f_i=d_Zh_i+h_i d_{F_iV}
$$
(and $d_{F_iV}=d_V|_{F_iV}$).
So we only need to construct $h_i$'s. This is done by induction, using Lemma
1.2.
Since $F_0V=0$, $f_0=0$ and we can take $h_0=0$. Now 1.2. guarantees that
having $h_i$ we can extend it to $h_{i+1}$. Namely, we first apply \cite{P2},
2.6.2, to the
semi split short exact sequence
$$
0@>>>F_iV@>>>F_{i+1}V@>>>Gr_{i+1}V@>>>0,
$$
so we can identify $F_{i+1}V$ with the cone over a morphism from $Gr_{i+1}V$ into $F_iV$. Now notice that the condition of
1.2.(ii) is met since $Gr_{i+1}V$ is K-projective and hence any morphism from
 $Gr_{i+1}V$ into $Z$ is homotopic to 0. So we can extend $h_i$ to a homotopy from $f_{i+1}$ to 0 and this extension is $h_{i+1}$.
\qed\enddemo

As explained above, this finishes the proof of the following result.

\proclaim{1.4. Theorem} Assume that $K$ is reductive. Then any 
$(\A,K,\D)$-module $V$ has a K-projective resolution $P@>>>V$ in $K(\A,K,\D)$. 
\qed\endproclaim

If we in addition assume that $\D$ is nonpositively graded
and that $V$ is bounded above by degree $k$, then in the above construction
each $P_{-i}$ can be taken bounded above by degree $k$. Namely, we can take $Q$
to be the classical projective resolution of $V$, i.e., $Q$ is bounded above
by $k$ with all $Q^i$ projective weak $(\A,K)$-modules (see A1.4.6).
Then $P_0$ is again bounded above by $k$ since $\D$ is nonpositively graded, 
etc. Now if all $P_{-i}$ are bounded above by degree $k$, then so is $P$.
In other words,

\proclaim{1.5. Corollary} Assume $\D$ is nonpositively graded.
Then any bounded above $(\A,K,\D)$-module $V$ has a K-projective resolution 
bounded above by the same degree as $V$.
\endproclaim

One can apply a dual construction to the one explained above to get existence
of enough K-injectives. Since this was done in \cite{BL2} (in essentially
the same way), we will not present this construction here.
Let us just comment on the difference between the two constructions.
The dual analogue of 1.2 is proved in the same way. For the dual analogue
of 1.3 we however need an additional finiteness assumption. Namely, the following is true: 

\proclaim{1.3'. Prop}
Let $V$ be an $(\A,K,\D)$-module with a decreasing Hausdorff filtration
$$
V=F_0V\supset F_1V\supset F_2V\supset\dots
$$
(Hausdorff means that $\cap_iF_iV=0$). Assume that
the graded objects $Gr_iV=F_iV/F_{i+1}V$
are K-injective for all $i$, 
that for any $i$, $F_iV$ is isomorphic to $F_{i+1}V\oplus
Gr_iV$ as a graded $(\A,K,\D)$-module,
and that for any $k\in\Bbb Z$, there are only finitely many $i$'s such that $Gr_i^kV\neq 0$.
Then $V$ is K-injective.
\endproclaim

The reason for this weakening of the result is the following. In 1.3, $V$ was 
the direct limit of
$F_iV$. Here, we want $V$ to be the inverse limit of $V/F_iV$, however this
is not true in general. More precisely, in the graded category, the inverse
limit of $V/F_iV$ is the algebraic part of the direct product 
$\prod_{k=0}^\infty Gr_kV$, 
while $V$ is the direct sum $\oplus_{k=0}^\infty Gr_kV$. These two are not
equal in general, but are of course equal if our finiteness assumption holds.

This difference now shows up in 1.1: to prove the dual claim we need a
finiteness assumption. The result is: if a complex $\Vd$ of
$(\A,K,\D)$-modules is bounded from below, if each $V^i$ is
K-injective, and if for any $k\in\Bbb Z$ the number of $i$'s such that
$(V^i)^{k-i}\neq 0$ is finite, then $s(\Vd)$ is K-injective.
Therefore we do not get an analogue of 1.4, but rather of 1.5.
Namely, we have to assume from the start that $\D$ is non-positively
graded, and that $V$ is bounded from below. Then $V$ embeds
quasiisomorphically into a K-injective $(\A,K,\D)$-module $I$, bounded
from below by the same degree as $V$.

\head 2. Bernstein-Lunts equivalence
\endhead

One of the main results of \cite{BL2}, Theorem 1.10, implies 
an equivalence of bounded equivariant
derived category $D^b(\A,K,N(\k))$ and the ordinary bounded derived category 
of Harish-Chandra modules $D^b(\M(\A,K))$, in case $\A$ is a projective
$\U(\k)$-module for the left multiplication composed with the map $\psi$. 
Their proof is rather complicated, since it uses K-injectives which are
complicated by themselves, and also in the proof they often need to dualize
various arguments. On the other hand, in \cite{BL1} there is an analogous result
10.12.5.1 for DG modules, which is rather simple, uses K-projectives, and
does not require boundedness. 
Of course, to have K-projectives the group $K$ has to be reductive. But
it turns out that once the result is proved for reductive groups, it is
easy to extend it to non-reductive case using the results of \cite{MP},
Section 2.

So let us first assume that $K$ is reductive. Then we know from Section 1 
that there are enough K-projective $(\A,K,\D)$-modules.
In our proof, it will be important to know that the forgetful functor
from the category of $(\A,K,\D)$-modules into the category of DG modules
over $\D$ preserves K-projectives. We start by establishing this
fact. The method to prove this is familiar: we show that this functor
has a right adjoint, which is `acyclic', i.e. maps acyclic complexes into
acyclic complexes (see Section A1.2).

We produce this functor in two steps. The first step is analogous to the
functor $\Ind_w$ from \cite{P2}, \S 1.3 and \cite{MP}, \S 2. Let $V$ be 
a DG module over $\D$ with action $\omega_V$ and differential $d_V$. We set
$\Ind_w(V)=R(K)\otimes V = R(K,V)$, where $R(K)$ is the algebra of regular
functions on $K$, and $R(K,V)$ is the space of regular functions from $K$ into
$V$. $K$ acts on $\Ind_w(V)$ by the right regular action, while
$x\in\D$ acts by 
$$
(\omega(x)F)(k)=\omega_V(\chi(k)x)(F(k)),\quad k\in K.
$$
The differential of $\Ind_w(V)$ is $1\otimes d_V$.
It is readily checked that in this way $\Ind_w(V)$ becomes a `weak 
$(K,\D)$-module', i.e., an algebraic $K$-module with a $K$-equivariant DG 
action of $\D$. Such modules can be identified with $(\U(\k),K,\D)$-modules,
if we let $\k\subset\U(\k)$ act by the difference of the two actions coming
from $K$ and $\D$.

\proclaim{2.1. Lemma} The above described functor 
$\Ind_w:\M(\D)@>>>\M(\U(\k),K,\D)$ is right adjoint to the forgetful functor.
\endproclaim

The proof of this is completely analogous to the proof that 
$\Ind_w:\M(\A)@>>>\M(\A,K)_w$ is right adjoint to the forgetful functor
(see \cite{MP}, 2.2). 

The second step consists of changing the algebra $\U(\k)$ to $\A$; here
the map $\U(\k)@>>>\A$ is given by the structural map $\psi:\k@>>>\A$. 
This is analogous to the functor $\pro_{\A,\B}$ mentioned in \cite{P2}, \S 1.2.
If $\gamma:\B@>>>\A$ is a $K$-equivariant morphism of algebras, we have a
forgetful functor from $\M(\A,K,\D)$ into $\M(\B,K,\D)$ given by the restriction
of scalars. Its right adjoint is given by
$$
V\longmapsto \pro_{\A,\B}(V)= \Hom_\B(\A,V)^{alg}.
$$
Here the $\B$-homomorphisms are with respect to the left multiplication on 
$\A$ (composed with $\gamma$). $\A$ acts on $\Hom_\B(\A,V)$ by the action $\pi$
which is right translation of the argument. $K$ acts by conjugation:
$$
\nu(k)f=\nu_V(k)\circ f\circ \phi_\A(k)^{-1},\quad k\in K.
$$
$x\in\D$ acts by $\omega(x)f=\omega_V(x)\circ f$. The differential is given
by $df=d_V\circ f$. Here, as usual, $\nu_V$, $\omega_V$ and $d_V$ are 
respectively the $K$-action, the $\D$-action, and the differential on $V$.
Finally, the algebraic part of $\Hom_\B(\A,V)$ is taken with respect to
the $K$-action; this is clearly invariant under all the above actions and the
differential. It is easy to check that $\pro_{\A,\B}(V)$ is an 
$(\A,K,\D)$-module, and that the functor $\pro_{\A,\B}$ is right adjoint to
the forgetful functor. The adjunction morphisms are just the standard ones
for extension of scalars. Namely, for an $(\A,K,\D)$-module $V$, 
$\Phi_V:V@>>>\pro_{\A,\B}(V)$ is given by
$$
\Phi_V(v)(a)=\pi_V(a)v, \quad a\in\A,\, v\in V,
$$
while for a $(\B,K,\D)$-module $W$, 
$\Psi_W:\pro_{\A,\B}(W)@>>> W$ is given by evaluation at $1$.

The composition of $\pro_{\A,\U(\k)}$ and $\Ind_w$ gives the desired right
adjoint to the forgetful functor from $\M(\A,K,\D)$ to $\M(\D)$. 
Furthermore, the same is true on the level of homotopic categories, as is
easy to check using \cite{P2}, \S 2.4.

It is obvious that $\Ind_w$ preserves acyclicity. We want to prove the same
for $\pro_{\A,\U(\k)}$. Here we need the assumption that $K$ is reductive,
but also an assumption on $\A$, namely that $\A$ is freely generated over
$\U(\k)$ by a (basis of a) $K$-submodule. This is a stronger assumption than 
the one in \cite{BL2}, where it is only assumed that $\A$ is projective over 
$\U(\k)$. However, our assumption is satisfied 
in the most interesting example: $\A=\U(\g)$, where $(\g,K)$ is a classical 
Harish-Chandra pair with $K$ reductive, i.e., $K$  
acts algebraically on $\g$ via inner automorphisms, $\k$ embeds into $\g$, 
and the differential of the $K$-action is $\ad$ composed with the  
embedding of $\k$ into $\g$; see \cite{P2}, Section 1.1.

\proclaim{2.2. Lemma} Let $(\g,K)$ be a Harish-Chandra pair, with $K$ reductive.
Then, as a $\U(\k)$-module for the left multiplication,
$$
\U(\g)=\U(\k)\otimes_\Bbb C \Cal P,
$$
where $\Cal P$ is a $K$-invariant subspace of $\U(\g)$.
\endproclaim
\pf
Let $\frak p$ be a $K$-invariant complement of $\k$ in $\g$. Let 
$X_1,\dots,X_n$ be a basis of $\k$ and let $Y_1,\dots Y_m$ be a basis of 
$\frak p$. Then
$$
X^iY^j;\quad i\in\Bbb Z_+^n,\,j\in\Bbb Z_+^m
$$
is a Poincar\'e-Birkhoff-Witt basis of $\U(\g)$. Let 
$\sigma:\Sl(\frak p)@>>>\U(\g)$ be the symmetrization map:
$$
\sigma(Z_1\dots Z_k)={1\over k!}\sum_{s\in S_k}Z_{s(1)}\dots Z_{s(k)}
$$
for $Z_1,\dots Z_k \in\frak p$. Then it is easy to check (see \cite{LM}, 
Lemma 2.2) that
$$
X^i\sigma(Y^j);\quad i\in\Bbb Z_+^n,\,j\in\Bbb Z_+^m
$$
is still a basis of $\U(\g)$, and that the span $\Cal P$ of $\{\sigma(Y^j);\,j\in\Bbb Z_+^m\}$ is $K$-invariant.
\qed\enddemo

Now if $\A=\U(\k)\otimes_\Bbb C \Cal P$ as a $\U(\k)$-module for the left 
multiplication, where $\Cal P$ is a $K$-invariant subspace of $\A$, then
as a $K$-module,
$$
\pro_{\A,\B}(V)=\Hom_\Bbb C(\Cal P,V)^{alg},
$$
depends only on the $K$-module structure of $V$. Since acyclicity of a complex
can be checked on the level of $K$-modules, it is enough to show that the
functor $V\longmapsto \Hom_\Bbb C(\Cal P,V)^{alg}$ from $\M(K)$ into $\M(K)$
is exact (here $\M(K)$ is the category of algebraic $K$-modules). This is
however true since $\M(K)$ is a semisimple category for reductive $K$.
Hence we have proved

\proclaim{2.3. Proposition} Assume that $K$ is reductive and that $\A$ is 
freely generated by a $K$-submodule as a $\U(\k)$-module for the left 
multiplication. Then the forgetful functor from the category of 
$(\A,K,\D)$-modules to the category of DG modules over $\D$ preserves 
K-projectives.
\endproclaim

We are ready now to prove the Bernstein-Lunts equivalence.  Let
$\eps:\D@>>>\E$ be a K-equivariant morphism of DG algebras, which is a
quasiisomorphism, i.e., induces an isomorphism in cohomology. We also
assume that the structural maps $\rho_\E:\k@>>>\E$ and
$\rho_\D:\k@>>>\D$ satisfy $\rho_E=\eps\circ\rho_D$.  The main example
is the counit map $N(\k)@>>>\Bbb C$, where $N(\k)$ is the standard
complex of the Lie algebra $\k$ (see \cite{P2}, 2.2.2 and 2.2.4).  Then 
$\eps$ induces
a forgetful functor $\For:\M(\A,K,\E)@>>>\M(\A,K,\D)$, which is just the
restriction of scalars. This functor has a left adjoint
$$
V\longmapsto V\otimes_\D \E
$$
as is proved in \cite{P2}, \S 2.5. Recall that the $\A$-action on 
$V\otimes_\D \E$ is the given action on the first factor, the $K$-action is 
on both factors, while the $\E$-action is given by the right multiplication
in the second factor, twisted to a left action. Both functors make sense
and remain adjoint on the level of homotopic categories.

Clearly, $\For$ preserves acyclicity, i.e., preserves
quasiisomorphisms, and therefore defines a functor on the level of
derived categories.  By 1.4, there are enough K-projective
$(\A,K,\D)$-modules, so by A1.4.4 and A1.3.3  the functor 
$V\mapsto V\otimes_\D \E$ 
has a left derived functor, which we denote by 
$V\mapsto V\dotimes_\D \E$. By A2.3.2, this functor is left adjoint to 
$\For:D(\A,K,\E)@>>>D(\A,K,\D)$.

\proclaim{2.4. Theorem} Assume that $K$ is reductive and that $\A$ is 
freely generated by a $K$-submodule as a $\U(\k)$-module for the left 
multiplication. Then the functors $\For$ and $V\mapsto V\dotimes_\D \E$, 
induced by $\eps:\D@>>>\E$ as above, are mutually inverse equivalences of 
categories $D(\A,K,\E)$ and $D(\A,K,\D)$.
\endproclaim
\pf
Let $V$ be an $(\A,K,\D)$-module, and let $P@>\delta>>V$ be a K-projective
resolution of $V$. Then $V\dotimes_\D \E=P\otimes_\D \E$, and the adjunction
morphism $\Phi_V:V@>>>\For(V\dotimes_\D \E)$ is given by the triple
$$
V {\underset ^{^{^{\sssize\sim}}} \to {\overset \,\,\,\delta\,\,\,\to \leftarrow}}
 P@>{id\otimes 1}>>\For(P\otimes_{\D}\E)
$$
($\sim$ denotes a quasiisomorphism). Namely, by A2.1.2, 
$\Phi_V=\alpha_{V,V\dotimes_\D\E}(1_{V\dotimes_\D\E})$; however, by the proof of 
A2.3.2, this can be identified with 
$$
V@<\delta<< P @>{\bar\alpha_{P,P\otimes_\D\E}(1_{P\otimes_\D\E})}>>
\For(P\otimes_\D\E),
$$
and $\bar\alpha_{P,P\otimes_\D\E}(1_{P\otimes_\D\E})=\bar\Phi_P= id\otimes 1$. 
Here $\bar\alpha$ and $\bar\Phi$ refer to the adjunction of $-\otimes_\D\E$ and 
$\For$.

So we only need to show that 
$id\otimes 1:P@>>>P\otimes_{\D}\E$ is a quasiisomorphism of 
$(\A,K,\D)$-modules. 
This morphism factors as
$$
P@>{id\otimes 1}>>P\otimes_{\D}\D@>{id\otimes\eps}>>P\otimes_{\D}\E,
$$
since $\eps(1)=1$. The first of these two morphisms is an isomorphism; its
inverse is given by the action map (this is the trivial change of DG algebras,
from $\D$ to $\D$). To show that the second morphism is a quasiisomorphism,
we can pass to DG modules over $\D$ (the property of being a quasiisomorphism 
can be checked on the level of complexes of vector spaces).

By 2.3, $P$ is a K-projective DG module over $\D$. Therefore
it is also a K-flat DG module over $\D$ (see \cite{BL1}, 10.12.4.4). 
In other words, the functor 
$V\mapsto P\otimes_\D V$, $V\in \M(\D)$ preserves acyclicity, or equivalently,
preserves quasiisomorphisms. Since $\eps:\D@>>>\E$ is a quasiisomorphism of
DG $\D$-modules, so is $id\otimes\eps:P\otimes_\D \D @>>> P\otimes_\D \E$.

So the adjunction morphism $\Phi_V$ is an isomorphism, for any 
$V\in D(\A,K,\D)$. Consider now the other adjunction morphism, 
$\Psi_W:\For W\dotimes_\D \E @>>> W$, for an $(\A,K,\E)$-module $W$.
The functor $\For$ clearly has the following property: a morphism $f$
in $D(\A,K,\E)$ is an isomorphism if and only if $\For f$ is an 
isomorphism in $D(\A,K,\D)$ (in other words, $\For$ preserves cohomology). 
Thus, to check that $\Psi_W$ is an isomorphism in the derived category,
it is enough to check that $\For\Psi_W$ is an isomorphism. However, by
adjunction, $\For\Psi_W\circ \Phi_{\For W} = id_{\For W}$. By the first part
of the proof, $\Phi_{\For W}$ is an isomorphism, hence so is $\For\Psi_W$.
This finishes the proof.
\qed\enddemo

We now specialize to the case mentioned before, namely the counit morphism
$\eps:N(\k)@>>>\Bbb C$. This is well known to be a quasiisomorphism, and
it clearly also satisfies our other assumptions. The category $D(\A,K,N(\k))$
is the equivariant derived category $D(\A,K)$ of $(\A,K)$-modules, and
$D(\A,K,\Bbb C)$ is the ordinary derived category $D(\M(\A,K))$ of 
$(\A,K)$-modules. The forgetful functor is the functor 
$Q:D(\M(\A,K))@>>>D(\A,K)$ which assigns to each complex of $(\A,K)$-modules
the same complex with all $i_\xi$, $\xi\in\k$ equal to $0$. So we get:

\proclaim{2.5. Corollary} Assume that $K$ is reductive and that $\A$ is 
freely generated by a $K$-submodule as a $\U(\k)$-module for the left 
multiplication. Then the forgetful functor $Q:D(\M(\A,K))@>>>D(\A,K)$ is 
an equivalence of categories.
\endproclaim

In this case we can eliminate the assumption of $K$ being reductive
(and also weaken the assumptions on $\A$), in the same way as in 
\cite{MP}, 2.14. Namely, it was proved in \cite{MP}, 2.11, that
the equivariant derived category $D(\A,K)$ is equivalent to the category
$D_{\M(\A,K)}(\A,L)$, where $L$ is a Levi factor of $K$. Here $D(\A,L)$ is the
equivariant derived category of $(\A,L)$-modules, and $D_{\M(\A,K)}(\A,L)$
is the full triangulated subcategory of equivariant $(\A,L)$-complexes with 
cohomology in $\M(\A,K)$. This was stated there for $\A=\U(\g)$, and 
$\A=\U_\theta$, the quotient of $\U(\g)$ corresponding to an infinitesimal 
character, but the proof clearly works for any $\A$. Namely, the inverse of 
the natural forgetful functor was given by the equivariant Zuckerman functor 
$R\Gamma^{equi}_{K,L}$, whose definition is independent of $\A$.

Assuming that $\A$ is freely generated over $\U(\l)$ by an $L$-submodule,
we can apply 2.5 to see that $Q_L:D(\M(\A,L))@>>>D(\A,L)$ is an equivalence
of categories. The restriction of $Q_L$ is then clearly an equivalence of
the category $D_{\M(\A,K)}(\M(\A,L))$ (i.e., the subcategory of $D(\M(\A,L))$ 
of complexes of $(\A,L)$-modules with cohomology in $\M(\A,K)$) with 
$D_{\M(\A,K)}(\A,L)$. Now $D_{\M(\A,K)}(\A,L)$ is equivalent to $D(\A,K)$
as is explained above. Analogously, $D_{\M(\A,K)}(\M(\A,L))$ is
equivalent to $D(\M(\A,K))$. This was proved in \cite{MP}, 1.12. It was 
stated there only for $\A=\U(\g)$, but the only assumption needed was that 
$\A$ is a flat $\U(\l)$-module for the right multiplication. In this case,
the proof of the Duflo-Vergne formula (\cite{MP}, 1.6) for 
$R^p\Gamma_{K,L}$ goes through without changes, and so then does the rest of 
\cite{MP}, \S 1.

Finally, it is obvious that the functor from $D(\M(\A,K))$ to $D(\A,K)$
induced by $Q_L$ via the above equivalences is precisely $Q_K$. So we get:

\proclaim{2.6. Corollary} Let $(\A,K)$ be a Harish-Chandra pair and $L$ a 
Levi factor of $K$. Assume that $\A$ is freely generated over $\U(\l)$ by 
an $L$-submodule. Then the functor $Q:D(\M(\A,K))@>>>D(\A,K)$ is an equivalence
of categories.
\endproclaim

\head A. Appendix: some facts from homological algebra
\endhead

The purpose of this appendix is to present some known things from homological
algebra needed in this paper, and some also in \cite{P2}, which are either 
unavailable in the literature, or exist but not in the best form for our purposes.
Some of the proofs are omitted or only sketched here, but they are written in 
more detail in \cite{P1}.
Most of them should in fact be easy to do, knowing the statements and looking
at analogous proofs in the "classical" situation, which is well covered by
the literature. 

In Section A1 we present a generalization of the well known construction of
derived functors between derived categories using adapted subcategories, to
the setting of triangulated categories (which are not necessarily homotopic 
categories of abelian categories), and their localizations (which are not
necessarily derived categories of abelian categories). 
In Section A2 we collect a few facts about
adjoint functors which are used over and over, both in this paper, and also
in \cite{P2}.

\head \bf A1. Derived Functors in Triangulated Category Setting \endhead

In this section, $\C$ will be a fixed triangulated category. For the
definition and basic properties of triangulated categories see \cite{Ve},
\cite{KS} or \cite{GM}. The main examples are homotopic and derived categories
of abelian categories; we will refer to these as the "classical situation".
The above books mostly treat this classical situation. Since equivariant
derived categories are not a priori derived categories of abelian categories, 
the results from the classical situation do not apply directly, and hence our 
need for a more general setting.

\subhead A1.1. S-systems and null systems \endsubhead
The following lemma contains several elementary facts that 
shall be needed below, but they are not explicitly stated in \cite{KS} or 
\cite{GM} except partly in the exercises. The proof is omitted here, but it 
is written out in detail in \cite{P1}.

\proclaim{A1.1.1. Lemma} Let $\C$ be a triangulated category. Then:

{\rm (i)} Suppose $X@>>>Y@>>>Z@>>>T(X)$ is a distinguished triangle in
$\C$. Then $X@>>>Y@>>>U@>>>T(X)$ is a distinguished triangle in
$\C$ if and only if $U$ is isomorphic to $Z$.

{\rm(ii)} Let $X@>f>>Y@>>>Z@>>>T(X)$ be a distinguished triangle in
$\C$. Then $f$ is an isomorphism if and only if $Z$ is isomorphic to
$0$.  

{\rm (iii)} Let $X_i@>f_i>>Y_i@>g_i>>Z_i@>h_i>>T(X_i)$, $i=1,2$, be two
distinguished
triangles in $\C$. Then their direct sum,
$X_1\oplus X_2@>f>>Y_1\oplus Y_2@>g>>Z_1\oplus Z_2@>h>>T(X_1\oplus X_2)$, is
also a distinguished triangle. Here $f=f_1\oplus f_2$, etc.

Conversely, if the direct sum is distinguished,
then both triangles are distinguished.

{\rm (iv)} For any two objects $X$ and $Y$ of $\C$, the triangle 
$$
\CD   X  @>0>>  Y  @>{i_2}>>T(X)\oplus Y  @>{p_1}>>  T(X)   
\endCD
$$
is distinguished.
\qed\endproclaim

An {\it S-system} in $\C$ is a
saturated localizing class $\Sl$ compatible with triangulation. The
saturation condition is 
$$ 
s\in\Sl \text{ if and only if } \exists \text{ morphisms } u,v 
\text{ such that }u\circ s\in\Sl ,\,\, s\circ v\in\Sl
$$
or, equivalently,
$$
s\in\Sl \text{ if and only if } Q(s) \text{ is an isomorphism. }
$$
Here $Q:\C@>>>\C[\Sl^{-1}]$ is the localization functor (recall that
$\C[\Sl^{-1}]$ is a triangulated category and $Q$ is exact).

This saturation condition is usually omitted. However, there are some technical
advantages of assuming it and on the other hand it is not restrictive. Namely, 
if $\Sl$ is any localizing class in $\C$ and $Q_\Sl:\C@>>>\C[\Sl^{-1}]$ the 
natural functor, one can define a saturated localizing class
$$
\overline\Sl=\{t\in \Hom\C\big| Q_\Sl(t)\text{ is an isomorphism }\}
$$
and $\C[\Sl^{-1}]$ and $\C[\overline\Sl^{-1}]$ are isomorphic.

A {\it null system} in $\C$ is a subfamily $\N$ of $\Ob(\C)$ satisfying the
following conditions:
\roster
\item"(N1)" $0\in\N$;
\item"(N2)" $X\in\N$ if and only if $T(X)\in\N$;
\item"(N3)" If $X@>>>Y@>>>Z@>>>T(X)$ is a distinguished triangle in
$\C$
and $X\in\N,Y\in\N$, then $Z\in\N$;
\item"(N4)" If $X\oplus Y\in\N$ then $X,Y\in\N$.
\endroster
In other words, the full subcategory generated by $\N$ is a
triangulated subcategory closed under isomorphisms
 and containing all direct summands of all its objects.
Namely the converse of (N4) holds because of (N2), (N3) and Lemma A1.1.1.(iv), 
so $\N$ is an additive subcategory. It is a triangulated subcategory by (N2) and
(N3). From (N3) and A1.1.1.(i) it follows that $\N$ is closed under
isomorphisms.

This definition differs from the one in \cite{KS}, which
does not require (N4) to hold. As we shall see, (N4) corresponds to
the saturation condition, so what we call a null system here could
be called a saturated null system.

\example{A1.1.2. Examples}
\roster
\item"(s1)" All isomorphisms in $\C$ form the smallest S-system in $\C$.
\item"(s2)" If $F:\C'@>>>\C$ is an exact functor and $\Sl$ an S-system in
$\C$,
then $F^{-1}(\Sl)=\{f\in\Hom\C'\big|F(f)\in\Sl\}$ is an S-system in
$\C'$.
\item"(s3)" All quasiisomorphisms in the
homotopic category of complexes over an abelian category form an S-system.
\item"(n1)" All objects of $\C$ isomorphic to $0$ form the smallest null system 
in $\C$, denoted by $\C_0$. 
\item"(n2)" If $F:\C'@>>>\C$ is an exact functor and $\N$ a null system in
$\C$,
then $F^{-1}(\N)=\{X\in\C'\big|F(X)\in\N\}$ is a null system in
$\C'$.
\item"(n3)" All acyclic complexes in the
homotopic category of complexes over an abelian category form a null system.
\endroster
\endexample

\proclaim{A1.1.3. Lemma} Let $\N$ be a null system in $\C$. Then for a
morphism $X@>f>>Y$ the following are equivalent:
\roster
\item"(i)" There is a distinguished triangle 
$X@>f>>Y@>>>Z@>>>T(X)$  with $Z\in\N$;
\item"(ii)" For any distinguished triangle 
$X@>f>>Y@>>>Z@>>>T(X)$, $Z$ is in $\N$.
\endroster
\endproclaim
\demo{Proof} Obvious from Lemma A1.1.1.(i) since $\N$ is closed under
isomorphisms.
\qed\enddemo

Given $\N$, let us denote by $\Sl(\N)$ the class of all morphisms of $\C$
satisfying the conditions of A1.1.3.  

\proclaim{A1.1.4. Proposition} $\Sl(\N)$ is an S-system. \endproclaim
\demo{Proof} 
It is proved in \cite{KS}, Prop.1.6.7., that $\Sl(\N)$ is
a localizing class. Compatibility with triangulation follows easily
from (N2) and (N3).
To show that $\Sl(\N)$ is also saturated, we need the following fact:

\proclaim{A1.1.5. Lemma} Let $Q:\C@>>>\C[\Sl(\N)^{-1}]$ be the localization 
functor. Then an object $X$ is in $\N$  if and only if  $Q(X)$ is isomorphic to 0.
\endproclaim
\demo{Proof}
Suppose that $Q(X)\cong 0$. It means that $1_{Q(X)}=0_{Q(X)}$, i.e., 
$Q(1_X)=Q(0_X)$. By an
elementary property of localization, this implies that there is a
morphism $Y@>s>>X$ in $\Sl(\N)$ such that $1_X\circ s=0_X\circ s$, 
that is $s=0$.
So $Y@>0>>X$ is in $\Sl(\N)$. Hence there is a distinguished triangle
$Y@>0>>X@>>>Z@>>>T(X)$ with $Z\in\N$. Using Lemma A1.1.1.(iv) and (i),
and the fact that $\N$ is closed under isomorphisms, we conclude that
$T(Y)\oplus X$ is in $\N$. However, then $X$ is also in $\N$ by (N4).

Conversely, let $X$ be in $\N$. Then from the distinguished triangle
$X@>1>>X@>>>0@>>>T(X)$ we conclude that the morphism $X@>>>0$ is in 
$\Sl(\N)$, so the corresponding morphism $Q(X)@>>>Q(0)=0$ is an
isomorphism.
\qed\enddemo
Now we can show that $\Sl(\N)$ is saturated. Suppose $Q(s)$ is an 
isomorphism. Let
$X@>s>>Y@>>>Z@>>>T(X)$
be a distinguished triangle. We have to prove that $Z$ is in $\N$.
However, applying $Q$ to the above triangle we get $Q(Z)\cong 0$ (by
A1.1.1.(ii)). So $Z\in\N$ by A1.1.5. This finishes the proof of A1.1.4.
\qed\enddemo

The following lemma is an easy consequence of A1.1.1.(ii).

\proclaim{A1.1.6. Lemma} Let $\Sl$ be an S-system in $\C$ and let $Q$ be the 
corresponding localization functor. Then the following are equivalent for an 
object $N$ of $\C$:
\roster
\item"(i)" $Q(N)\cong 0$;
\item"(ii)" For any distinguished triangle $X@>f>>Y@>>>N@>>>T(X)$, 
$f$ is in $\Sl$;
\item"(iii)" There is a distinguished triangle $X@>f>>Y@>>>N@>>>T(X)$
with $f\in\Sl$.
\qed\endroster\endproclaim

Given $\Sl$, let $\N(\Sl)$ be the full subcategory of all objects satisfying the
conditions of A1.1.6. We shall sometimes call these objects $\Sl$-acyclic.

\proclaim{A1.1.7. Proposition} $\N(\Sl)$ is a null system. 
\endproclaim
\demo{Proof} Follows from A1.1.2.(n1), (n2). Namely, $\N(\Sl)$ is the 
inverse under $Q$ of the null system 
$\C[\Sl^{-1}]_0$.
\qed\enddemo

\proclaim{A1.1.8. Theorem} Attaching $\Sl(\N)$ to $\N$ and $\N(\Sl)$ to
$\Sl$ gives a one-to-one correspondence between S-systems and null
systems in $\C$. 
\endproclaim
\demo{Proof}
Clearly $\N\subset\N(\Sl(\N))$ and $\Sl\subset\Sl(\N(\Sl))$.
$\N(\Sl(\N))\subset\N$ is Lemma A1.1.3. So it remains to prove that
$\Sl(\N(\Sl))\subset\Sl$. However, $s\in\Sl(\N(\Sl))$ implies the existence
of a distinguished triangle $X@>s>>Y@>>>N@>>>T(X)$ with $N\in\N(\Sl)$.
So $Q_\Sl(N)\cong 0$ and therefore $Q_\Sl(s)$ is an isomorphism by A1.1.1.(ii).
Hence $s\in\Sl$ since $\Sl$ is saturated.
\qed\enddemo

\example{A1.1.9. Remark} 
Under the correspondence from A1.1.8., 
the examples from A1.1.2. correspond to each other as follows:
\roster
\item A1.1.2.(s1) corresponds to A1.1.2.(n1) (by Lemma A1.1.1.(ii).)
\item A1.1.2.(s2) corresponds to A1.1.2.(n2), meaning that
if $\N$ and $\Sl$ in $\C$ correspond to each other, then $F^{-1}(\N)$ and 
$F^{-1}(\Sl)$ also correspond to each other. 
\item A1.1.2.(s3) corresponds to A1.1.2.(n3) .
\endroster
\endexample

We also remark that instead of null systems one can consider {\it thick 
subcategories} of $\C$, as it is done in \cite{Ve} and \cite{GM}.
It is however easy to prove that this notion is the same as our notion of
a null system.

\subhead A1.2. Derived functors and acyclic functors \endsubhead
Let $\C$ be a triangulated category, $\Sl$ an S-system in $\C$, $\N$
the corresponding null system of $\Sl$-acyclic objects and
$Q:\C@>>>\C[\Sl^{-1}]$ the localization functor.

Let $(\C',\Sl',\N')$ be another triple as above, and let $F:\C@>>>\C'$
be an exact functor. We call $F$ $(\Sl,\Sl')$-acyclic if $F(\Sl)\subset 
F(\Sl')$, or equivalently $F(\N)\subset F(\N')$. These two properties
are equivalent by the description of the
correspondence between S-systems and null systems.

As an example, consider two abelian categories $\A$ and $\A'$. Let
$\Sl$ and $\Sl'$ be the S-systems of all quasiisomorphisms in $K(\A)$
and $K(\A')$, and let $\N$, $\N'$ be the null systems of acyclic
complexes (see A1.1.9.(3)). Then for any functor $F:\A@>>>\A'$ which
is exact (in the classical sense), $K(F)$ is $(\Sl,\Sl')$-acyclic.

For the rest of this section we assume that $\N'$ is the null system $\C'_0$ of
all objects isomorphic to 0 and $\Sl'$ the corresponding S-system of
all isomorphisms (see A1.1.9.(1)). In this case, F is called just 
$\Sl$-acyclic. 

Let $F:\C@>>>\C'$ be an arbitrary exact functor.
A right derived functor of $F$ is a pair
$(RF,\varepsilon_F)$, where
$$
RF:\C[\Sl^{-1}]@>>>\C'
$$
is an exact functor and
$$
\varepsilon_F:F@>>>RF\circ Q
$$
is a morphism of functors\footnote{The definition of a functor $F$ between triangulated 
categories includes a fixed isomorphism from $TF$ into $FT$, where $T$ is the translation
functor. A morphism of functors is assumed to be compatible with these isomorphisms.}, 
such that the following universal property
holds: if $(G,\varepsilon)$ is another such pair, then there exists a
unique
morphism of functors $\eta:RF@>>>G$ such that the diagram
$$
\CD   F  @>{\varepsilon_F}>> RF\circ Q        \\
    @V=VV              @V{\eta\circ Q}VV    \\
      F  @>{\varepsilon}>>   G\circ Q  
\endCD
$$
commutes. If $RF$ exists, it is unique up to isomorphism. This 
follows from the universal property in a standard way.

Dually, one defines a left derived functor of $F$:  it is a pair
$(LF,\varepsilon_F)$, where $LF:\C[\Sl^{-1}]@>>>\C'$ 
is an exact functor and $\varepsilon_F:LF\circ Q @>>>F$
is a morphism of functors, such that the following universal property
holds: if $(G,\varepsilon)$ is another such pair, then there exists a
unique
morphism of functors $\eta:G@>>>LF$ such that 
$\varepsilon=\varepsilon_F\circ (\eta\circ Q)$.

These definitions in this
generality are due to P.Deligne, \cite{De}.

A right (or left) derived functor of $F$ does not have to exist. 
If however $F$ is $\Sl$-acyclic, then by the universal property of localization
there is a unique functor $\bar{F}:\C[\Sl^{-1}] @>>> \C'$
such that $F=\bar{F}\circ Q$. Furthermore, $\bar{F}$ is exact,
and it is easy to prove that

\proclaim {A1.2.1. Proposition} If $F$ is $\Sl$-acyclic, then
$\bar{F}$
is a right derived functor of $F$, if we take
$\varepsilon_F=1:F@>>>\bar{F}\circ Q$ (and also left derived, with
$\varepsilon_F=1:\bar{F}\circ Q@>>>F$).
\qed\endproclaim

Since $\bar F$ is just the factorization of $F$ through the
localized category, it is usually
denoted again by $F$ to simplify notation. 

\example{A1.2.2. Remark} In the more general situation when we
have an arbitrary S-system $\Sl'$ in $\C'$, we consider the composition
$Q_{\Sl'}\circ F$ and we define the right (left) derived functor of
$F$ to be the right (left) derived functor of $Q_{\Sl'}\circ F$. In
the latter case the above definition applies.
\endexample

\subhead A1.3. Adapted subcategories \endsubhead
If $F$ is not acyclic, one tries to find a triangulated subcategory
$\D$ of $\C$ such that $F\big|_{\tsize\D}$ is acyclic. Then $F$ will
factor through the localization of $\D$ with respect to the induced S-system 
$\Sl_\D$. In case $\D$ is big enough, every object of $\C[\Sl^{-1}]$ will be
isomorphic to an object of $\D[\Sl_\D^{-1}]$ and we will be able to
construct $RF$ and $LF$ using that. This is the same construction as in the 
classical situation, but we will review it in some detail to show that it 
works in our more general situation.

If $\D$ is a triangulated subcategory of $\C$
then $\D$ is a triangulated category with the induced
structure: a triangle in $\D$ is distinguished if and only if it is distinguished
as a triangle in $\C$ and the translation functor on $\D$ is the restriction
of the translation on $\C$.

\proclaim {A1.3.1. Lemma} Let $\D$ be a triangulated subcategory of $\C$,
$\Sl$ an S-system in $\C$ and $\N$ the corresponding null system. 
Then $\N_\D=\N\cap\Ob\D$ is a null system in $\D$ and the
corresponding S-system $\Sl_\D$ is equal to $\Sl\cap\Hom \D$.
\endproclaim
\demo {Proof} This is just a special case of A1.1.2.(n2) and A1.1.9.(2)
where the functor $F$ is the inclusion from $\D$ into $\C$.
\qed\enddemo

Now we can form the localized category $\D[\Sl_\D^{-1}]$ and we get an
exact functor $\Psi: \D[\Sl_\D^{-1}]@>>>  \C[\Sl^{-1}]$. Namely, the elements
of $\Sl_\D$ clearly go to isomorphisms under $\D@>>>\C@>>>\C[\Sl^{-1}]$. 

Let us now further assume that $\D$ satisfies one of the following
properties:
\roster
\item"($b_R$)" For any $X\in\C$ there exists $X@>s>>M$ with $s\in\Sl$ and
$M\in\D$
\item"($b_L$)" For any $X\in\C$ there exists $M@>s>>X$ with $s\in\Sl$ and
$M\in\D$
\endroster

\proclaim {A1.3.2. Lemma} Let $\D$ be a triangulated subcategory of $\C$
such
that ($b_R$) or ($b_L$) holds. Then:
\roster
\item"(i)" $\Psi$ is fully faithful, so $\D[\Sl_\D^{-1}]$ can be
viewed
as a full subcategory of $\C[\Sl^{-1}]$;
\item"(ii)" $\D[\Sl_\D^{-1}]$ is a triangulated subcategory of $\C[\Sl^{-1}]$
and its natural triangulated structure (coming from $\D$ and
localization) is equal to the induced triangulated structure;
\item"(iii)" $\Psi$ is an equivalence of categories and any of its
quasiinverses is exact.
\endroster
\endproclaim
\demo{Proof} Done in the classical situation.
\qed\enddemo

\proclaim{A1.3.3. Theorem} Let $\D$ be a triangulated subcategory of $\C$
satisfying ($b_R$) or ($b_L$) above. Let $F:\C@>>>\C'$ be an exact
functor such that 
\roster
\item"(a)" The restriction of $F$ to $\D$ is $\Sl_\D$-acyclic.
\endroster
Then $F$ has a left derived functor $(LF,\varepsilon_F)$ in case ($b_L$) holds
(and $F$ has a right derived functor in case ($b_R$) holds). Furthermore, for any
$M\in\D$, $\varepsilon_F(M)$ is an isomorphism.
\endproclaim
\demo{Proof} Analogous to the classical situation.
\qed\enddemo

We shall say that $\D$ is left (right) adapted to $F$ if it 
satisfies the conditions (a) and ($b_L$) ((a) and ($b_R$)).

\subhead A1.4. $\Sl$-projective and $\Sl$-injective objects \endsubhead

Let $\C$ be a triangulated category, $\Sl$ an S-system in $\C$
and $\N$ the corresponding null system of $\Sl$-acyclic objects.

\proclaim{A1.4.1. Theorem} Let $P$ be an object of $\C$. Then the
following conditions are equivalent:
\roster
\item For any $X\in\N$, $\Hom_\C(P,X)=0$;
\item If $s:X@>>>P$ is in $\Sl$, then there exists a
$t\in\Hom_\C(P,X)$
such that $s\circ t=1_P$ (we will also see that $t$ is unique and in $\Sl$);
\item For any $X\in\Ob\C$, the natural homomorphism between
abelian groups $\Hom_\C(P,X)$ and $\Hom_{\C[\Sl^{-1}]}(P,X)$
mapping $P@>f>>X$ into $P@<1<<P@>f>>X$ is an isomorphism;
\item For any diagram
$$
\CD   
    @.  X \\
   @.  @VVsV \\
   P @>f>> Y
\endCD
$$
in $\C$ with $s\in\Sl$, there is a unique $g\in\Hom_\C(P,X)$ such that
$f=s\circ g$.
\endroster
\endproclaim
\demo{Proof}  See \cite{Sp}, 1.4. and the remark below.
\qed\enddemo

\example{A1.4.2. Remark} Suppose there is an exact bifunctor $H$ from
$\C^{opp}\times\C$ into the homotopic category of complexes of abelian
groups, such that for any two objects $X,Y$ in $\C$
$$
H^i(H(X,Y))=\Hom_\C(X,Y[i]).
$$
For example in the classical situation we can take $H=\Hom^{\cdot}$.
Then clearly the condition (1) in A1.4.1. is equivalent to the condition
\roster
\item"(1')" The functor $H(P,-)$ maps acyclic objects into acyclic
complexes of abelian groups (in the classical sense).
\endroster
This is the condition that is used in \cite{Sp}.
\endexample

Objects that satisfy the conditions of A1.4.1. are called $\Sl$-projective
objects. Dually, one defines $\Sl$-injective objects. They were first
introduced by Verdier in \cite{Ve}; he calls them free on the left,
respectively right. Their basic properties and some applications were 
studied by Spaltenstein in \cite{Sp}. He considered the classical case of
the homotopic category $K(\A)$ of an abelian category $\A$, with $\Sl$
being the class of all quasiisomorphisms and $\N$ the null system of
acyclic complexes. To suggest that the definitions are intrinsic to the 
homotopic category, he called  $\Sl$-projective objects of $K(\A)$ K-projective,
and $\Sl$-injective objects K-injective. The same terminology is used by
Bernstein and Lunts in \cite{BL1} and \cite{BL2} in the case of homotopic
category of equivariant complexes. We will follow this terminology in all these 
cases, but in the present generality the prefix `K' does not make sense. 

The importance of $\Sl$-projectives and $\Sl$-injectives lies in the fact that
they can be used as adapted categories to define derived functors, as
we are going to explain now.  

\proclaim{A1.4.3. Proposition} The full subcategory $\Cal P$ of $\C$ consisting
of $\Sl$-projective objects is a null system in $\C$. The same holds
for $\Cal I$, the full subcategory of $\Sl$-injectives. In particular,
$\Cal P$ and $\Cal I$ are triangulated subcategories of $\C$.
\endproclaim
\demo{Proof} This follows easily from property (1) in A1.4.1. Namely $\Cal P$
is additive since the functor $\Hom_\C(-,X)$ is additive. $\Cal P$ is closed under
$T$ and $T^{-1}$ since the null system of $\Sl$-acyclic objects is closed under 
$T$ and $T^{-1}$, and furthermore $\Hom_\C(TP,X)=\Hom_\C(P,T^{-1}X)$ and 
$\Hom_\C(T^{-1}P,X)=\Hom_\C(P,TX)$.

Since
the functor $\Hom_\C(-,X)$ is cohomological, whenever
$P@>>>Q@>>>K@>>>T(P)$ is a distinguished triangle in $\C$ and
$P$ and $Q$ are in $\Cal P$, then $K$ is also in $\Cal P$. Finally,
using the additivity of $\Hom_\C(-,X)$ again, we see that 
a direct summand of an $\Sl$-projective object is 
$\Sl$-projective.

The proof for $\Sl$-injectives is analogous.
\qed\enddemo

\proclaim{A1.4.4. Theorem} Suppose that there are enough $\Sl$-projectives
in $\C$, i.e., that $\Cal P$ satisfies the condition $(b_L)$ from 
Section A1.3. Then $\Cal P$ is left adapted to any exact functor
$F$ from $\C$ to another triangulated category $\C'$. An analogous
claim is true for $\Sl$-injectives.
\endproclaim
\demo{Proof}  
It is enough to check the condition (a) from Theorem A1.3.3. However, if
$P$ is $\Sl$-projective and $\Sl$-acyclic, then $\Hom_\C(P,P)=0$, which
implies that $P$ is isomorphic to 0 in $\C$. Hence $F(P)$
must also be isomorphic to 0.
\qed\enddemo 

In the classical situation of homotopic categories over abelian categories,
one usually uses resolutions by complexes with
injective or projective components. Let us clarify how this fits into the
framework of the theory studied here. As announced earlier, in this situation
we call $\Sl$-injective objects K-injective, and $\Sl$-projective objects
K-projective.

\proclaim{A1.4.5. Proposition} Let $\A$ be an abelian category. Let $P^\cdot$
be a complex over $\A$, bounded above, such that all $P^j$ are projective 
objects of $\A$. Then $P^\cdot$ is a K-projective object of $K^-(\A)$ and
$K(\A)$.

Dually, if $I^\cdot$ is a complex over $\A$, bounded below and having injective 
components, then $I^\cdot$ is a K-injective object of $K^+(\A)$ and $K(\A)$.
\endproclaim
\pf To show that $P^\cdot$ is a K-projective object of $K^-(\A)$, we have to 
show that if $\Cd$ is an acyclic complex bounded above, then
any chain map from  $P^\cdot$ to $\Cd$ is homotopic to zero. However, this is a
well-known fact from classical homological algebra, used to prove that
left derived functors are well defined. 

To get the claim for $K(\A)$, we can just note that since $P^\cdot$ is 
bounded above, the chain maps from $P^\cdot$ into any complex $\Cd$ are actually
chain maps into an appropriate truncation of $\Cd$. However, if $\Cd$ is acyclic,
so is any of its truncations. 

The other claim is proved in the same way.
\qed\enddemo

An example in the introduction of \cite{Sp} (taken from \cite{Do}) shows that 
A1.4.5 is not true without the boundedness assumptions. 
The following fact is well-known:

\proclaim{A1.4.6. Theorem} Let $\A$ be an abelian category. If $\A$ has enough
projectives, then for any $\Xd$ in $K^-(\A)$ there is a complex $P^\cdot$ with
projective components, bounded above by the same degree as $\Xd$, and mapping
quasiisomorphically onto $\Xd$. In particular, $K^-(\A)$ has enough K-projectives.

Dually, if $\A$ has enough
injectives, any $\Xd$ in $K^+(\A)$ can be quasiisomorphically embedded into
a complex with injective components, bounded below by the same degree
as $\Xd$. In particular, $K^+(\A)$ has enough K-injectives.
\qed\endproclaim

It is shown in \cite{Sp} that, under certain 
conditions, there are also enough K-injectives and K-projectives in $K(\A)$.
In particular, we need

\proclaim{A1.4.7. Theorem} {\rm (\cite{Sp}, 3.5)}
If an abelian category $\A$ has enough projectives, has direct limits, and
the direct limit functor is exact, then $K(\A)$ has enough K-projectives.
\qed\endproclaim

\head A2. Some Properties of Adjoint Functors \endhead

In this section we review some standard and some not so standard facts
about adjoint functors. Most of the general theory can be found for example 
in \cite{ML}. The case of abelian categories is easy and well known. I do not
however know any reference for A2.3.2; I learned it from D.~Mili\v ci\'c.

\subhead A2.1. Definition and general properties \endsubhead
Let $\A$ and $\B$ be two categories. Let $F$ be a functor from $\A$ to $\B$
and let $G$ be a functor from $\B$ to $\A$. We say that $F$ is left adjoint
to $G$, or that $G$ is right adjoint to $F$, if the bifuctors $\Hom_\B(F(-),-)$
and $\Hom_\A(-,G(-))$ from $\A^{opp}\times\B$ into $\Cal Sets$ are isomorphic.
This condition means that
for any two objects $X$ from
$\A$ and $Y$ from $\B$, there is a bijection
$$
\alpha=\alpha_{X,Y}:\Hom_\B(FX,Y)@>\cong>>\Hom_\A(X,GY),
$$
natural in $X$ and $Y$. It is easy to show that if $F$ and $F'$ are both left 
adjoint to a functor $G$ from $\B$ to $\A$, then $F$ and $F'$ are isomorphic. 
In other words, we can speak of the (unique) left adjoint of $G$.
Analogously, right adjoints are also unique.

Another easy fact is about the adjoint of a composition. Namely,
suppose $F:\A@>>>\B$ and $F':\B@>>>\C$ are two functors, with right adjoints 
$G$ and $G'$ respectively. Then $GG'$ is right adjoint to $F'F$. 

A related fact with obvious proof is the following:

\proclaim{A2.1.1. Proposition}
Assume $I:\A@>>>\B$ is fully faithful, i.e., $\A$ can be identified with a full
subcategory of $\B$.
Let $\C$ be another category, and let $F:\C@>>>\A$ be a functor. 
If $H$ is left adjoint to $IF$, then $HI$ is left adjoint to $F$.
Similarly, if $H$ is right adjoint to $IF$, then $HI$ is right adjoint to $F$.

Note that $HI$ can be viewed as the restriction of $H$ to $\A$. 
\qed\endproclaim

\proclaim{A2.1.2. Theorem} Let $F$ from $\A$ to $\B$ and $G$ from $\B$ to $\A$
be two functors. Then $F$ is left adjoint to $G$ if and only if there are
natural transformations
$$
\Phi:Id_\A@>>>GF,\quad \Psi:FG@>>>Id_\B,
$$
such that for any object $X$ of $\A$, the composition
$$
FX@>{F(\Phi_X)}>>FGF(X)@>{\Psi_{FX}}>>FX
$$
is the identity morphism, and for any object $Y$ of $\B$, the composition
$$
GY@>{\Phi_{GY}}>>GFG(Y)@>{G(\Psi_Y)}>>GY
$$
is the identity morphism.
\endproclaim
\pf If $\alpha:\Hom_\B(F(-),-)@>\cong>>\Hom_\A(-,G(-))$ gives adjunction of $F$
and $G$, we define $\Phi$ and $\Psi$ by
$\Phi_X=\alpha_{X,FX}(1_{FX})$
and $\Psi_Y=\alpha_{GY,Y}^{-1}(1_{GY})$.

Conversely, given $\Phi$ and $\Psi$, define $\alpha$ by 
$\alpha_{X,Y}(\phi)=G(\phi)\circ\Phi_X$. The inverse is given by
$\beta_{X,Y}(\psi)=\Psi_Y\circ F(\psi)$.

In filling in the details, naturality plays a key role.
\qed\enddemo

The morphisms $\Phi_X$, $X\in\A$ and $\Psi_Y$, $Y\in\B$ are called adjunction
morphisms. The following two statements are easy consequences of A2.1.2; recall
that a $\Bbb C$-category is an additive category such that the $\Hom$-groups are
not only abelian groups, but have
an additional structure of a vector space over $\Bbb C$ (compatible with the
group structure), and the 
composition law is bilinear. In the case of $\Bbb C$-categories, additive 
functors are required to be linear on morphisms.

\proclaim{A2.1.3. Corollary} Let $F$ from $\A$ to $\B$ be left adjoint to $G$ 
from $\B$ to $\A$. If $\A$ and $\B$ are additive categories and if $F$ or $G$ is
an additive functor, then the other functor is also additive, and the maps
$\alpha_{X,Y}$ are isomorphisms of abelian groups. 

An analogous claim is true for $\Bbb C$-categories.
\qed\endproclaim

\proclaim{A2.1.4. Proposition} Let $F:\A@>>>\B$ be left adjoint to $G$ and let 
$\Phi$ and $\Psi$ be the corresponding natural transformations. Then:
\roster
\item"(i)" $\Phi_X$ is an isomorphism for every $X\in\A$ if and only if
$F$ is fully faithful. In that case, $G$ is essentially onto.
\item"(ii)" $\Psi_Y$ is an isomorphism for every $Y\in\B$ if and only if
$G$ is fully faithful. In that case, $F$ is essentially onto.
\qed\endroster
\endproclaim

\example{A2.1.5. Examples} There are many examples of adjoint functors in
various branches of mathematics. Let us mention just a few: the functor 
assigning to a set $X$
the free group generated by $X$ is left adjoint to the forgetful functor
from groups to sets; analogous claims are true in case of algebras, modules, etc.
The functor assigning to a Lie algebra its universal enveloping algebra is
left adjoint to the forgetful functor from associative algebras into Lie 
algebras. Limits and colimits also serve as examples of adjunction; see
\cite{P1}, Section 2.2. 

In the following, let us discuss in more detail three classes of examples that
are of special interest for us.

\smallskip

\no (1) Equivalences of categories.

Let $F:\A@>>>\B$ be an equivalence of categories. Let
$G:\B@>>>\A$ be a quasiinverse for $F$, i.e., $GF$ and $FG$ are isomorphic to 
the identity functors. Then both $F$ and $G$ are fully faithful and 
essentially onto; in particular, we have natural isomorphisms
$$
\Hom_\B(FX,Y)\cong\Hom_\A(GFX,GY)\cong\Hom_\A(X,GY).
$$
So $F$ is left adjoint to $G$. Analogously, $F$ is right adjoint to $G$. It is
now clear from A2.1.4 that all the corresponding adjunction morphisms are 
isomorphisms.

Conversely, if $F$ is left adjoint to $G$ and if all adjunction morphisms
$\Phi_X$ and $\Psi_Y$ are isomorphisms, then $F$ and $G$ are clearly
mutually inverse equivalences of categories.

\smallskip

\no (2) Duality functors.

Let $\A$ be a category and $D:\A@>>>\A$ a contravariant functor. We say that
$D$ is a duality on $\A$, if for any $X$ and $Y$ in $\A$ there is an isomorphism
$$
\gamma_{X,Y}:\Hom_\A(X,DY)@>\cong>>\Hom_\A(Y,DX),
$$
natural in $X$ and $Y$. Note that $D$ defines two covariant functors, 
$D':\A^{opp}@>>>\A$ and $D'':\A@>>>\A^{opp}$. The above definition then says that
$D'$ is right adjoint to $D''$. If we adopt the convention which identifies
$D$ with $D'$, we can say that $D$ is right adjoint to its opposite functor.

The adjunction morphisms corresponding to this adjunction, when interpreted
in $\A$, both give the same morphism $\Phi_X:X@>>>DDX$. 

Familiar dualities usually satisfy much stronger conditions, like $DD$ being
isomorphic to the identity functor, and in case of abelian $\A$, $D$ being
exact. An example of a weaker duality is
the standard duality for infinite dimensional vector spaces. It is not hard
to show that the natural inclusion $X@>>>X^{**}$ gives adjunction as above,
but $X^{**}$ is not isomorphic to $X$. Another example is the duality for
finitely generated modules over a commutative ring (or algebra) $A$:
$$
D(M)=\Hom_A(M,A).
$$
This is not even an exact functor; however it becomes an equivalence if we
pass to the derived category.

\smallskip

\no (3) Restriction and extension of scalars.

This is another well-known example. Many constructions in this paper and also
in \cite{P2} are
direct generalizations of this example.

Let $f:A@>>>B$ be a morphism of algebras (or just rings). Then there is
an obvious `restriction of scalars' functor, or forgetful functor from the
category $\M(B)$ of (left) modules over $B$ into the category $\M(A)$: we 
simply let $a\in A$ act on a $B$-module $M$ as $f(a)$. This functor $\For$ has
both adjoints:

The left adjoint is the functor 
$$
M\longmapsto B\otimes_A M
$$
where $B$ is viewed as a right $A$-module via right multiplication composed
with $f$, and $B$ acts on $B\otimes_A M$ by left multiplication in the first
factor. On morphisms, the functor is given by $f\mapsto 1\otimes f$. 

The adjunction morphisms are given as follows: For $M\in\M(A)$, define
$\Phi_M:M@>>> \For(B\otimes_A M)$ by $\Phi_M(m)=1\otimes m$.
For $N\in\M(B)$, define 
$\Psi_N:B\otimes_A\For N@>>>N$ by $\Psi_N(b\otimes n)=bn$,
where $bn$ denotes the action of $b$ on $n$. One now checks that these are well
defined morphisms in appropriate categories, and that they
really give the required adjunction.

The right adjoint is given by
$$
M\longmapsto\Hom_A(B,M)
$$
for $M\in\M(A)$. Here $B$ is considered as a left $A$-module via left 
multiplication composed with $f$, and the action of $B$ on $\Hom_A(B,M)$ is 
given by right multiplication in the first variable.

The adjunction morphisms are defined similarly as for the left adjoint, using
the action map and the evaluation at 1.
\qed\endexample

To end this section, let us just mention the fact that if $F$ is left
adjoint to $G$, then $F$ preserves colimits (sums, cokernels, pushouts,
direct limits, etc.), while $G$ preserve limits (products, kernels, pullbacks,
inverse limits etc.).
An important special case of this is the case of functors between
abelian categories; then $F$ is right exact (since it preserves
cokernels), and $G$ is left exact (since it preserves kernels). These facts
are well known; a detailed proof together with the definitions and some more 
related facts can be found in \cite{P1}, \S 2.2.

\subhead A2.2. The case of abelian categories \endsubhead

In this section, $\A$ and $\B$ are abelian categories, and $F:\A@>>>\B$ is left
adjoint to $G:\B@>>>\A$. As we already mentioned at the end of last section,
we have

\proclaim{A2.2.1. Proposition} $F$ is right exact and $G$ is left exact.
\endproclaim

Recall that an object $X$ of $\A$ is called projective if the functor 
$\Hom_\A(X,-)$ is exact. $X$ is called injective if the functor $\Hom_\A(-,X)$ is
exact.
The following result is very simple, but it is crucial for many constructions
in homological algebra.

\proclaim{A2.2.2. Theorem} If $G$ is exact, $F$ preserves projectives. If $F$ is
exact, $G$ preserves injectives.
\endproclaim
\pf
Let $X\in\A$ be projective. Then
$$
\Hom_\B(FX,-)=\Hom_\A(X,G(-))
$$
is exact as a composition of two exact functors. Hence $FX$ is projective.
The other claim is proved in the same way.
\qed\enddemo

A2.2.2 gives a standard way of constructing projectives and injectives in abelian
categories. We say that an abelian category $\A$ has enough projectives if
every object of $\A$ is a quotient of a projective object. $\A$ has enough
injectives if every object of $\A$ embeds into an injective object.

\proclaim{A2.2.3. Theorem} Suppose that $G$ is exact, that $\A$ has enough
projectives, and that the adjunction morphism $\Psi_B$ is an epimorphism for 
every object $B$ of $\B$. Then $\B$ has enough projectives.

Dually, if $F$ is exact, $\B$ has enough injectives, and $\Phi_A$ is a 
monomorphism for every $A\in\A$, then $\A$ has enough injectives.
\endproclaim
\pf
Let $B\in\B$. Since $\A$ has enough projectives, there is an epimorphism
$P@>>>GB$ in $\A$, with $P$ projective. Applying $F$,
we get an epimorphism $FP@>>>FGB$, since $F$ is right exact by A2.2.1. Composing
this with $\Psi_B:FGB@>>>B$ we get an epimorphism $FP@>>> B$. However, $FP$ is 
projective by A2.2.2. Hence $\B$ has enough projectives. The second claim is 
proved in the same way.
\qed\enddemo

\subhead A2.3. The case of triangulated categories \endsubhead

Let us now go back to the situation of Section A1. Let $\C$ and $\D$ be two 
triangulated categories with null systems $\N$ and $\M$ and corresponding
S-systems $\Sl$ and $\T$. Let us denote by $\C_\Sl$ and $\D_\T$
the corresponding localizations.

Let $F:\C@>>>\D$ and $G:\D@>>>\C$ be exact functors. Assume that $F$ is left 
adjoint to $G$.

\proclaim{A2.3.1. Theorem} 
If $G$ is $(\T,\Sl)$-acyclic (i.e., sends $\M$ into $\N$), then 
$F$ maps $\Sl$-projective objects of $\C$ into $\T$-projective objects of $\D$.
Dually, if $F$ is $(\Sl,\T)$-acyclic then $G$ maps $\T$-injectives into 
$\Sl$-injectives.
\endproclaim
\pf Let $P\in\C$ be $\Sl$-projective and let $X\in\M$. Then
$$
\Hom_\D(F(P),X)=\Hom_\C(P,G(X))=0
$$ 
since $G(X)\in\N$. So $F(P)$ is $\T$-projective. The proof of the dual
statement is analogous.
\qed\enddemo

Next, we want to show that if $F$ and $G$ have derived functors $LF$ and $RG$,
then these derived functors are still adjoint. In fact, to prove this we
need to assume that $LF$ and $RG$ can be calculated using adapted subcategories. 
Of course, this will be true in all practical situations. The proof presented
here is due to D.~Mili\v ci\'c.

\proclaim{A2.3.2. Theorem} Let $F:\C@>>>\D$ be left adjoint to $G:\D@>>>\C$.
Let $\L\subset\C$ be a subcategory left adapted to $F$ and let $\R\subset\D$
be a subcategory right adapted to $G$. Let $LF:\C_\Sl@>>>\D_\T$ be the left 
derived functor of $F$, and let $RG:\D_\T@>>>\C_\Sl$ be the right derived 
functor of $G$; these functors exist under the above assumptions. Then $LF$ is 
left adjoint to $RG$.
\endproclaim
\pf Let $X\in\C_\Sl$ and $Y\in\D_\T$. We have to prove that
$$
\Hom_{\D_\T}(LF(X),Y)\cong\Hom_{\C_\Sl}(X,RG(Y)),
$$
naturally in $X$ and $Y$. By our assumptions, $X$ is naturally isomorphic
in $\C_\Sl$ to an object $A$ of $\L$, while $Y$ is naturally isomorphic in
$\D_\T$ to an object $B$ of $\R$. Also, $LF(X)=F(A)$ and $RG(Y)=G(B)$. 
Therefore it is enough to show that
$$
\Hom_{\D_\T}(F(A),B)\cong\Hom_{\C_\Sl}(A,G(B)),
$$
naturally in $A$ and $B$.

Let $\phi:F(A)@>>>B$ be a morphism in $\D_\T$. We can represent $\phi$ by a
triple
$$
F(A)@>f>>C@<s<< B
$$
where $f$ and $s$ are morphisms in $\D$, and $s\in\T$. 
We can assume $C\in\R$, passing to an equivalent triple if necessary.

Let $\gamma=\gamma_{A,C}$ be the map from $\Hom_\D(F(A),C)$ into 
$\Hom_\C(A,G(C))$ defined by the adjunction of $F$ and $G$. So $\gamma$ is 
natural in $A$ and $C$. Since $B$ and $C$ are in $\R$, $G(s)$ is in $\Sl$. 
Therefore the triple
$$
A@>{\gamma(f)}>>G(C)@<{G(s)}<< G(B)
$$
represents a morphism in $C_\Sl$ from $A$ to $G(B)$. We denote this morphism
by $\alpha(\phi)$. 
One shows $\alpha(\phi)$ is well-defined, i.e., that it
depends only on the morphism $\phi$ and not on the choice of $f$ and
$s$. 
Furthermore, $\phi\mapsto\alpha(\phi)$ is natural in $A$ and $B$.

Finally, in order to prove that 
$\alpha:\Hom_{\D_\T}(F(A),B)@>>>\Hom_{\C_\Sl}(A,G(B))$ is a bijection, we 
construct its inverse. Let $\psi:A@>>>G(B)$ be a morphism in $\C_\Sl$. We can 
represent it by
$$
A@<s<< C@>f>>G(B)
$$
with $s\in\Sl$. As before, we can assume $C\in\L$. We define $\beta(\psi)$ to
be the class of the triple
$$
F(A)@<{F(s)}<<F(C)@>{\gamma^{-1}(f)}>>B.
$$
One now shows that $\beta$ is indeed inverse to $\alpha$.
\qed\enddemo

We will mostly use the above theorem in special cases, when each of the 
functors $F$ and $G$ is either acyclic (so we can take the whole category as
an adapted category), or its derived functor can be calculated
using $\Sl$-projectives, respectively $\T$-injectives. In those cases one could 
produce slightly simpler proofs for the theorem, but this does not seem to
compensate for the loss of generality.

\example{A2.3.3. Remark} Let us show how to apply A2.3.2 to the classical 
situation.
Let $F:\A@>>>\B$ be additive and left adjoint to $G:\B@>>>\A$, where $\A$ and
$\B$ are abelian categories. Then it is a standard (and easy) fact that $F$
and $G$ define exact functors on the level of homotopic categories $K(A)$ and 
$K(B)$. It is immediate that $F$ is still left adjoint to $G$ on the homotopic
level. However, to define $LF$ and $RG$, one usually has to restrict $F$ to
$K^-(\A)$ and $G$ to $K^+(\B)$ to get the existence of adapted subcategories.
In that case A2.3.2 does not even make sense. If however
one of the functors, say $F$, is cohomologically bounded, then $LF$ can be
extended to $D(\A)$ and then restricted to $D^+(\A)$. Inspecting this 
construction one sees that $LF:D^+(\A)@>>>D^+(\B)$ is still calculated using a 
left adapted subcategory. Therefore A2.3.2 can be applied in this situation.

Another possibility is that $LF$ and $RG$ can be defined
on full derived categories because of existence of enough K-projectives in 
$K(\A)$ and K-injectives in $K(\B)$. Then A2.3.2 applies without any assumptions
on $F$ or $G$.
\endexample

\head Acknowledgements \endhead

Most of the material in this paper is taken from my 1995 University of
Utah Ph.D. thesis, written under guidance of my advisor Dragan 
Mili\v ci\'c. I would like to thank him for all the help, suggestions
and ideas he unselfishly shared with me.

\enddocument
\end